\UseRawInputEncoding 

\documentclass{article}

\usepackage{color,graphicx,times}
\usepackage{amssymb,amscd}
\usepackage{pb-diagram}     
\usepackage{amsmath}
\usepackage{amsthm}
\usepackage{graphicx}

\makeatother

\let \ttorg \tt \def \tt{\ttorg \obeyspaces}

\begin{document}

\title{\Large\bf Topology of Vortex Reconnection}
\author{Louis H. Kauffman\\ Department of Mathematics, Statistics \\ and Computer Science (m/c
249)    \\ 851 South Morgan Street   \\ University of Illinois at Chicago\\
Chicago, Illinois 60607-7045\\ $<$kauffman@uic.edu$>$\\}

\maketitle

\thispagestyle{empty}

\subsection*{\centering Abstract}

{\it This paper studies how knotted vortices undergo reconnection in their process of simplification. We consider the world line of such a knotted vortex and show how the relationships
of knot theory and four dimensional topology can be used to give lower bounds on its reconnection number and how to explicitly determine the reconnection number for torus links and positive links.}
\bigbreak

\noindent {\bf Keywords.} Alexander polynomial, Alexander-Conway polynomial,  Seifert pairing, Seifert surface, Seifert circles, Khovanov homology, knotted vortex, reconnection, reconnection number, knot signature, Milnor conjecture,Rasmussen invariant.

\section{Introduction}
This paper discusses knotted vortices and how they tend to degenerate to unknots and unlinks by reconnection processes. We define a reconnection number of knots and links and show how some aspects 
four dimensional topology allow the computation of the reconnection number. This allows us the see that in some cases of experiments with computer models of vortex knots, the models do find the least number of reconnections to undo the knot. In other cases the physical process is less efficient than the best mathematical reconnection. Our main result is an exact determination of the reconnection number $R(K)$ for $K$ any positive link.  In particular, we show:\\

{\it   If $K$ is an unsplittable knot or link with a positive diagram with $c(K)$ crossings and $s(K)$ Seifert circles, then the reconnection number of $K$ is exactly twice the genus of the Seifert spanning surface of $K$ when $K$ is a knot, and for both knots and links we have the formula $$R(K) = c(K) - s(K) +1.$$ }\\

 If $K = T(a,b)$ is a torus knot of type $(a,b)$ then the reconnection number of $K$ is equal to the product $(a-1)(b-1).$\\
 
The meanings of the above terms are given in the body of the paper.\\

This mathematical result about the topology of vortex reconnection is proved herein by using  deep results about four-ball surfaces bounding knots and links in three-space. These results are  due to Rasmussen \cite{Rasmussen} via his use of  Khovanov homology.
Results related to those of Rasmussen were originally proved by Kronheimer and Mrowka using gauge theory and so have their roots in relationships of topology and physics. Kronheimer and Mrowka proved a conjecture of John Milnor about the torus knots and singularities of algebraic varieties in a space of two complex variables. We return to physics in this paper by
applying a result about genera of surfaces in the four-ball to properties of vortices in the four dimensions of spacetime.  In this way we also return to the ideas that generated knot theory at the hands of Lord Kelvin
and his nineteenth century theory of atoms as vortices in the luminiferous aether.\\

The interest in these results about reconnection numbers for knots is both mathematical and physical. It is a very nice mathematical question to ask how many oriented saddle point reconnections are needed to 
unknot a knot. This is nearly as basic a question as the query about the unknotting number of knot -- how many crossing switches will unknot the knot. On the physical side, vortices in fluids and superfluid vortices tend to unknot via physical reconnection. Knowing the least possible numbers of reconnections needed to undo a knot makes possible better comparisons with the way such knots behave in
experiments.\\

\section{Vortex Knot Reconnection}
\begin{figure}
     \begin{center}
     \begin{tabular}{c}
     \includegraphics[width=8cm]{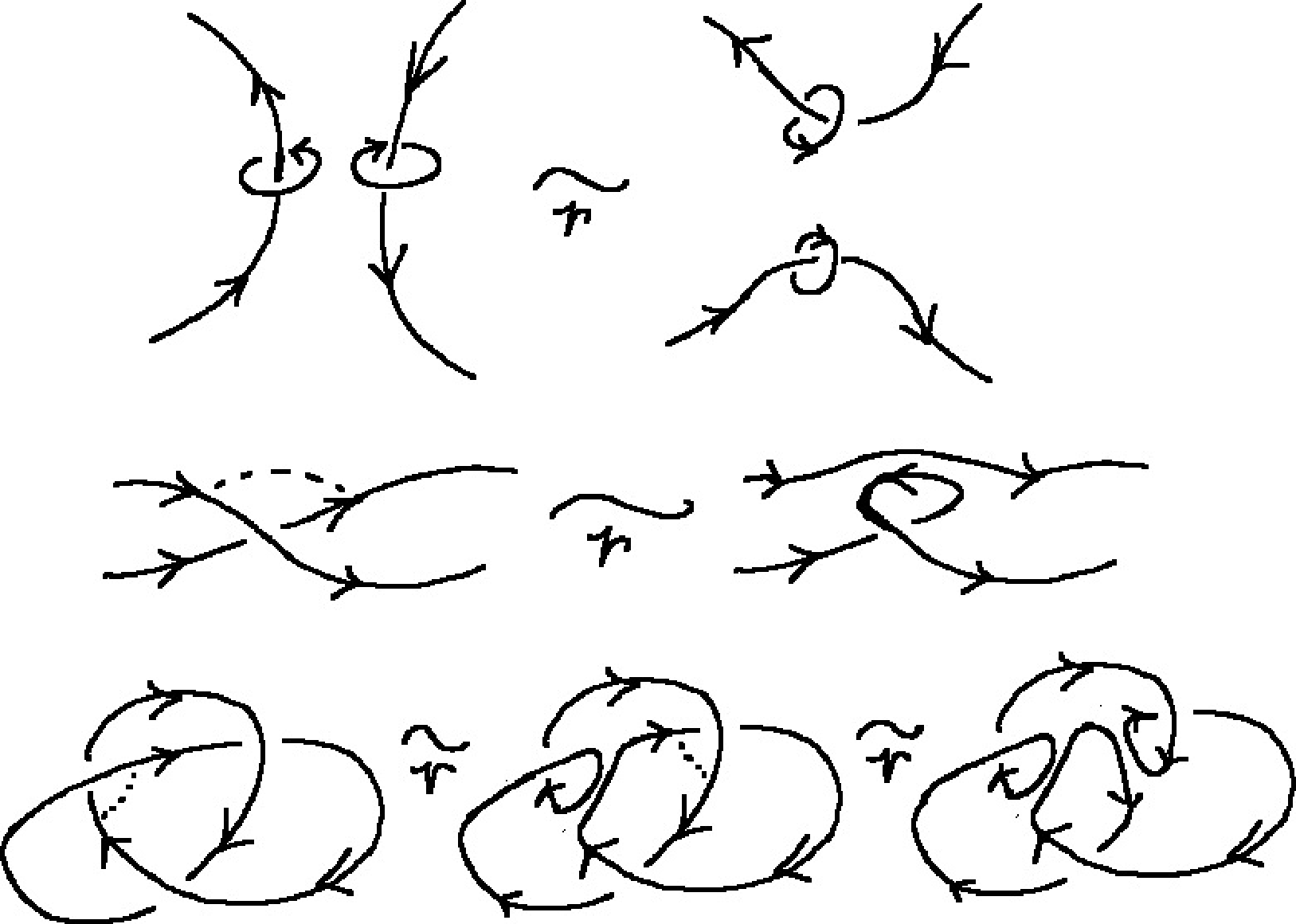}
     \end{tabular}
     \caption{\bf Diagrammatic Vortex Reconnection and Topology Change}
    \label{recon}
\end{center}
\end{figure}

\begin{figure}
     \begin{center}
     \begin{tabular}{c}
     \includegraphics[width=13cm]{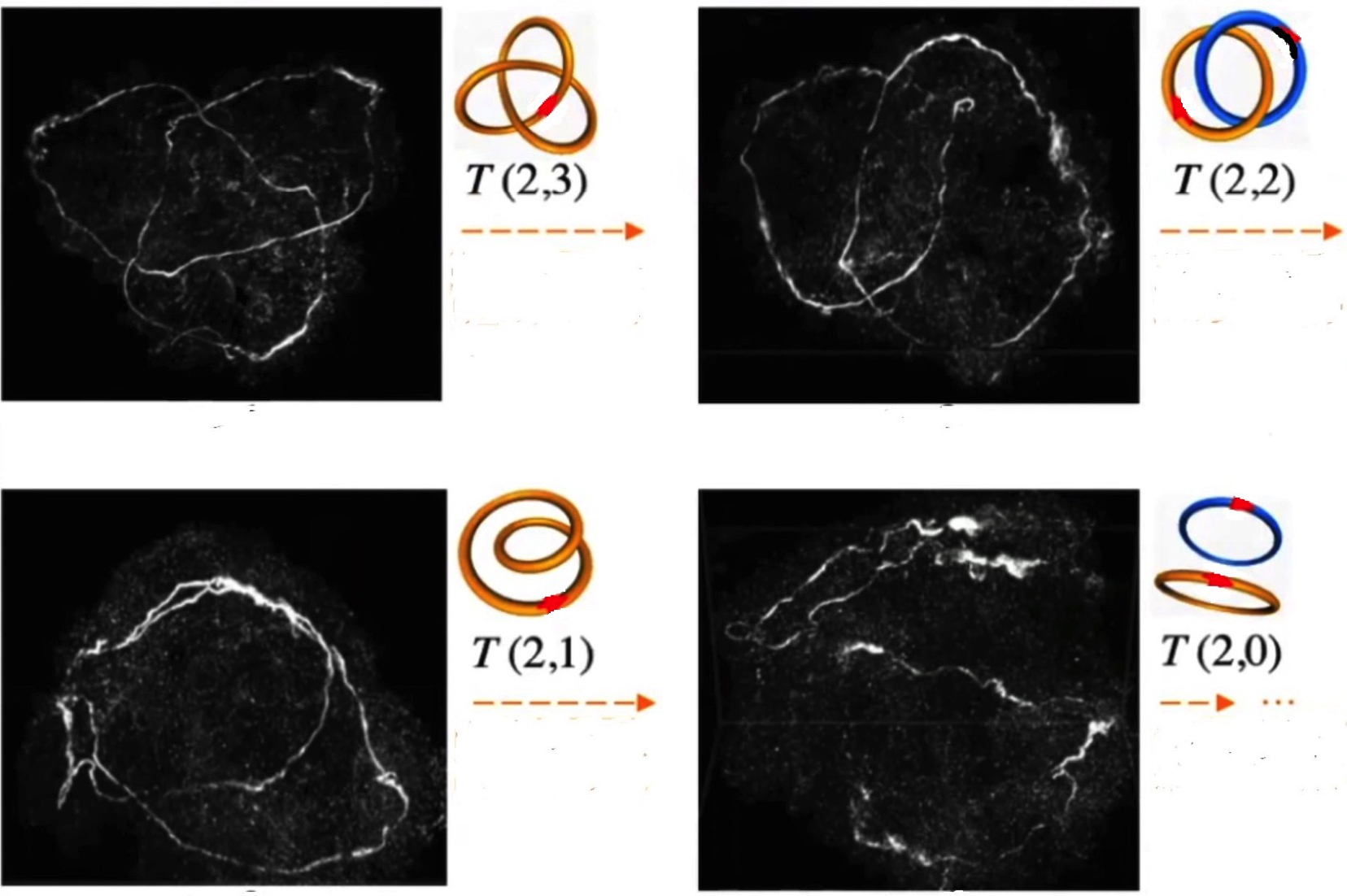}
     \end{tabular}
     \caption{\bf Trefoil cascade}
    \label{tcascade}
\end{center}
\end{figure}

\begin{figure}
     \begin{center}
     \begin{tabular}{c}
     \includegraphics[width=10cm]{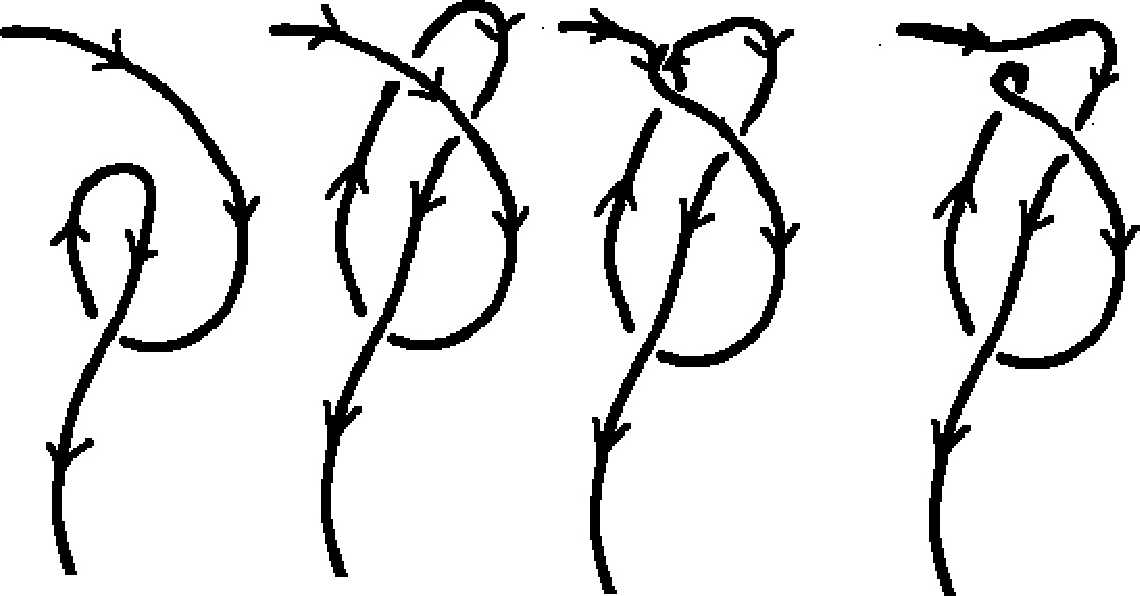}
     \end{tabular}
     \caption{\bf Vortex Reconnection Producing a Link}
    \label{vortexrecon}
\end{center}
\end{figure}

\begin{figure}
     \begin{center}
     \begin{tabular}{c}
     \includegraphics[width=13cm]{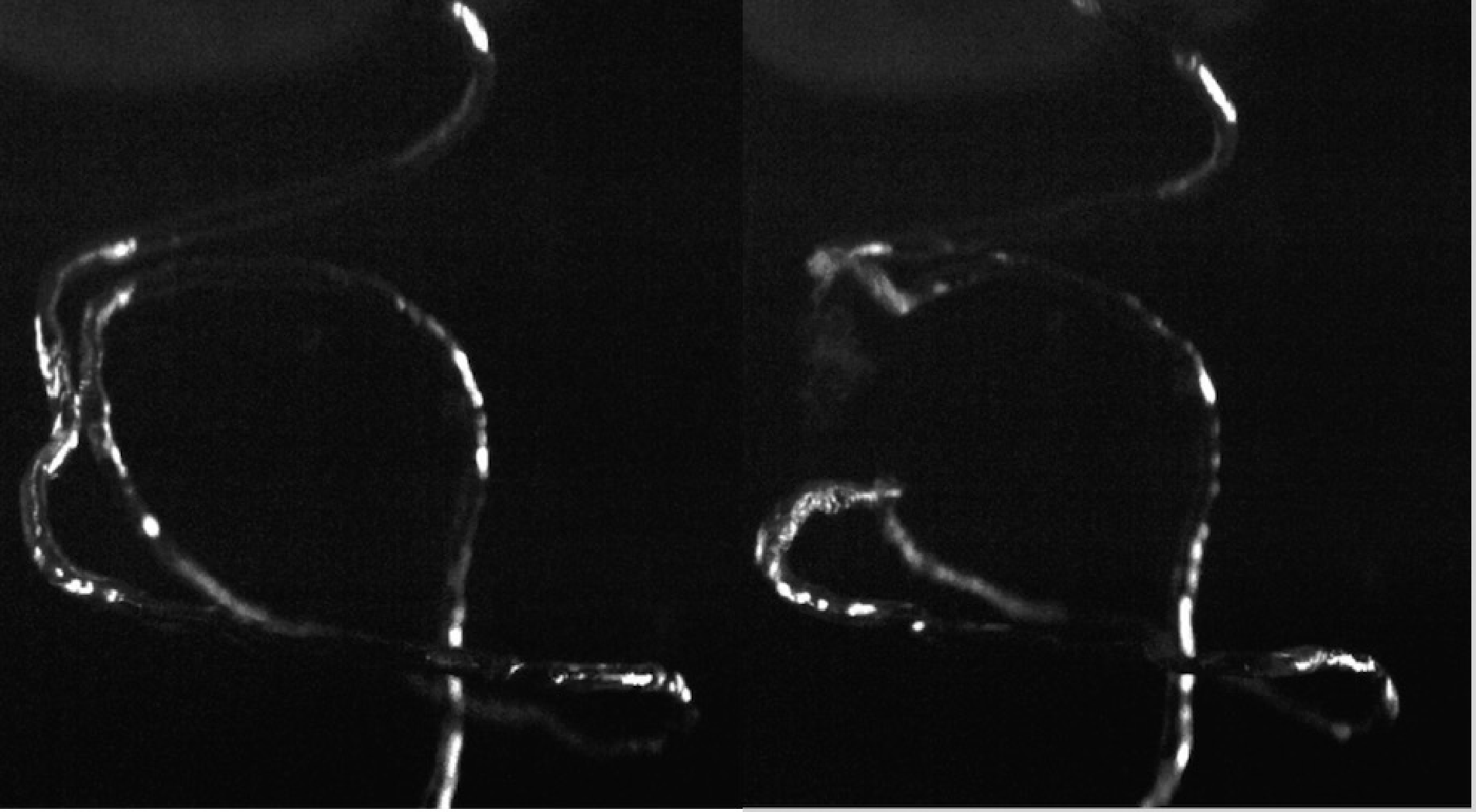}
     \end{tabular}
     \caption{\bf Vortex Reconnection and Topology Change}
    \label{turbine}
\end{center}
\end{figure}

We consider knotted vortices. The main result of the paper is an exact determination of the {\it reconnection number} $R(K)$ for $K$ any {\it positive} link. The meanings of these terms will be made clear in 
the body of the paper.  Reconnection  is a process that happens spontaneously to a one-dimensional vortex, creating a recombination of the vortex lines and a change in topology. The reconnection number is 
the least number of reconnections needed to transform the knot to an unknotted loop. Figure~\ref{recon} gives a diagrammatic illustration of reconnection. Figure~\ref{tcascade} and Figure~\ref{turbine} are photographs of actual
reconnection of vortices in water \cite{KI,Aleks}. The Figure~\ref{tcascade} illustrates the experimental work of Kleckner and Irvine \cite{KI} and shows a water vortex in the form of a trefoil knot undergoing a cascade of
reconnections that result in unknotted and unlinked circles in parallel with the diagrammatic illustration in Figure~\ref{recon}. Figure~\ref{vortexrecon} is a drawing of the vortex reconnection that occurs in the photograph
in Figure~\ref{turbine}.  In these two figures we see that a reconnection can, under appropriate circumstances, cause the topology to become more complicated. In this case a loop linking the main vortex is produced by the reconnection. This occurs in a turbine experiment performed by Professor Alexander Alekseenko \cite{Aleks}.
For other work on reconnection, see \cite{Ricca1,Ricca2,Ricca3}.\\

Lord Kelvin (Sir William Thompson) in the nineteenth century theorized that material atoms were knotted vortices in the luminiferous aether.
His theory led him to request the making of tables of knots, and such tables were constructed by Peter Guthrie Tait, T. P. Kirkman, and C. N. Little. These early tables preceded the development of knot theory in terms of topological invariants, and were a strong motivation for their development. The Kelvin theory of knotted vortex atoms  
did not survive the criticisms of the aether. These criticisms were driven by the negative results of the Michelson-Morley experiment and the advent of Einstein's theory of special relativity. The possibility of knotting at the level of elementary particles has not been ruled out. 
Knotting of fluid vortices remains of great interest to physicists. For many years the possibility of knotted vortices in three dimensional 
fluid dynamics was a topic of theoretical interest, but no one had explicitly exhibited knotted vortices in even such a ubiquitous fluid as water. This situation changed in 2012 with the work of William Irvine and 
Dustin Kleckner \cite{KI}. They produced repeatable experiments that yield knotted vortices in water. Their method is a generalization of earlier methods that produce ring vortices in smoke and indeed in water.
Their method is to make a knotted template and to pull it quickly through still water. The edge of the template produces a knotted vortex that can be observed by high speed photography. The resulting visual evidence of
knotted vortices in a palpable fluid is fascinating and very suggestive of many questions about the behaviour of such structures.\\

One can  simulate knotted vortices and their dynamics using the Gross-Pitaevskii non-linear Schroedinger equation \cite{KKI}. In both experiments with actual fluids and in computer experiements it is seen that the knotted vortices
change by a process called {\it reconnection}. At the diagram level 
a reconnection is what we have called a re-smoothing in our discussions of knot and link invariants.  Two parallel vortex filaments come close to one another. An exchange process
occurs that results in a transition to a new topological connection corresponding to a re-smoothing as illustrated in Figure~\ref{recon}.   The vortices must be locally spinning (around the one dimensional vortex line)  in the same direction for reconnection to occur. This directionality is illustrated in the Figure~\ref{recon}. In  this figure the reader will see a directionality along the vortex line and a spin of the vortex indicated by the right hand rule with respect to the vortex line. This same direction of spin is reinforced locally. Note that two parallel vortex lines have opposite directionality. In Figure~\ref{turbine} we illustrate photographs from a turbine experiment performed by Professor Alexander Alekseenko \cite{Aleks}.   In the first photograph the vortex line is connected. In the seond photograph the reconnection has divided the vortex into two components with one component encircling the other. The topology has made a transition to the linking of two vortex lines.\\

\begin{figure}
     \begin{center}
     \begin{tabular}{c}
     \includegraphics[width=7cm]{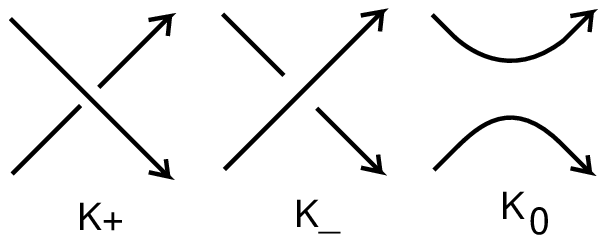}
     \end{tabular}
     \caption{\bf Crossing Signs in a Skein Triple}
     \label{skein}
\end{center}
\end{figure}

There are consequences of reconnection for the related combinatorial topology. In particular, if we consider oriented knot and link diagrams under regular isotopy (Reidemeister moves two and three)  and also oriented reconnection as described in the last paragraph, then one sees that the {\it writhe} of the diagram is invariant under these moves. Recall that for a diagram $D$, the writhe, $wr(D)$ is the sum of the crossing signs in the diagram. In Figure~\ref{skein} the signs of crossings are indicated, with the diagram $K_{+}$ having a sign of $+1$ and the diagram $K_{-}$ having a sign of $-1.$ In the Figure~\ref{recon}
we illustrate the result of reconnection near a diagrammatic crossing. The reader will have no difficulty seeing that the writhe is preserved under reconnection. At the bottom of the figure we illustrate two successive
reconnections on a trefoil diagram. After the first reconnection, the diagram has transformed into a link of two circles. One more reconnection transforms this link to a connected and unknotted loop.\\

How many reconnections are needed to unknot a knot, or to unlink a link? This is a natural mathematical question and it is a key physical question about knotted vortices since their observed behaviour is to undergo reconnection and become a collection of unlinked and unknotted circles in three dimensional space. See \cite{KKI} for experimental results about this question. In this section of the paper, we will construct lower bounds on this {\it reconnection number}  by regarding the process of reconnection as the creation of a saddle point in a surface that represents the world line of the evolving vortex. Then the reconnection number can be compared with the genus of a world line surface for the vortex. A nontrivial application of four dimensional topology gives the lower bound.\\ 

We say that two oriented knots are {\it cobordant} if they form the ends of an embedded tube ($S^{1} \times I$) in the cartesian product of three dimensional space with a unit interval. One can think of a cobordism as a process that transforms one knot into the other. The process consists in a sequence of births of unknotted circles, deaths of unknotted circles and passage through oriented saddle points. The saddle point passages correspond to local oriented arc reconnections as indicated by the figures. A process that creates a genus zero surface
(a topological tube) is called a {\it concordance}. We can consider surfaces in the four-ball $D^{4}$ that bound a knot embedded in the boundary three-sphere $S^{3}.$ If $F$ is a surface in $D^{4}$ that bounds a knot $K$ in $S^{3}$, then $F$ can be described by a process of births and deaths of circles and passage through saddle points. For a given process, we can calculate the genus of the surface that is produced. The {\it four ball genus $g_{4}(K)$} is the least genus among all surfaces that bound the knot $K$ in the four-ball. From the point of view of processes, the four-ball genus is the least genus traced by any process that creates a surface bounding the knot. There are general results known about the four-ball genus. For example, it is known \cite{Mura} that 
$| \sigma(K)|/2 \le  g_{4}(K)$ where $\sigma (K)$ denotes the classical signature of the knot. This is computed as the signature of $V+V^{T}$ where $V$ is a Seifert matrix for the knot and will be discussed in Section 3.
If $K$ is a link with $\mu$ components, then the result generalizes to 
$$| \sigma(K)| \le  2 g_{4}(K) + \mu - 1.$$
Also, if $u(K)$ is the unknotting number of a link $K$ (the least number of crossing switches needed to convert $K$ (from any diagram or in space) to the unlink, then
$$g_{4}(K) \le u(K).$$
\bigbreak

The relevance of the process view of cobordisms of knots is that the saddle passage is exactly a writhe-preserving oriented reconnection. In studying vortices one is interested in the least number of such reconnections that will transform a knot to a collection of unknotted circles. Let us define $R(K)$ to be this {\it reconnection number for the knot $K$}. That is $R(K)$ is the least number of reconnections that will transform $K$ to a collection of unknotted circles. By thinking in terms of the four-ball genus, we see that $g_{4}(K) \le R(K)/2,$ since pairs of saddle moves contribute to the increment of one handle, whence an addition of one to the genus. Thus we find the result\\

\noindent {\bf Theorem.} Let $\sigma(K)$ denote the signature of an oriented knot $K$ and $R(K)$ the reconnection number of $K.$ Then $| \sigma(K) | \le R(K).$ \\

The background and proof of this result will be discussed below. For example, if $K$ is the trefoil knot, then $|\sigma(K)| = 2,$
from which we conclude that it will take at least two reconnections to undo the trefoil, and indeed two reconnections suffice.
\bigbreak

Cobordisms as defined below correspond to allowed reconnections of vortices.
Writhe is preserved by these reconnections. 
The deaths and births of unknotted circles are relevant to the cobordism topology but not directly to the vortices.
Bounds on the number of saddle points that will make a cobordism from  a knot to a collection of unknots can be computed from the topology. We will detail more about the topology and how to use it to study reconnection in the next sections.\\

\begin{figure}
     \begin{center}
     \begin{tabular}{c}
     \includegraphics[width=8cm]{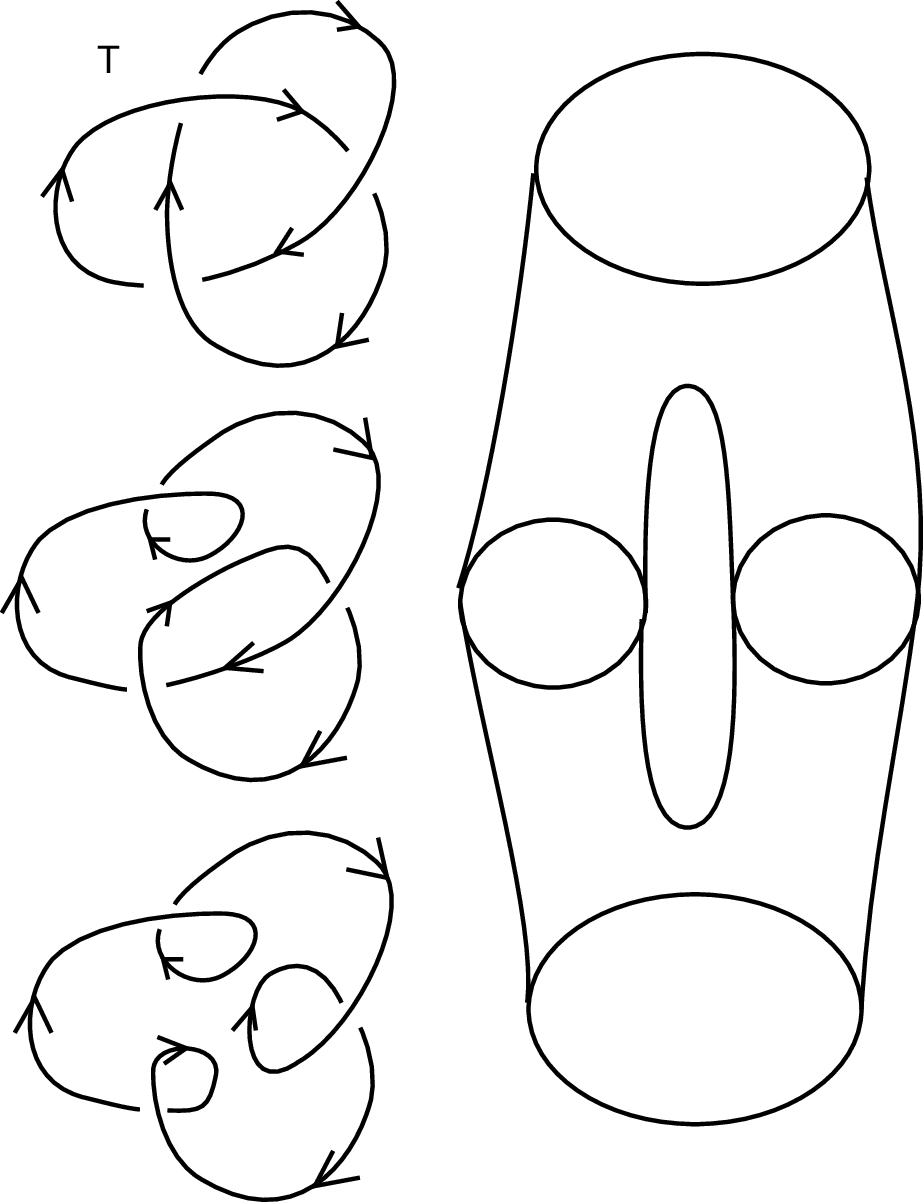}
     \end{tabular}
     \caption{\bf A Genus One Reconnection Sequence}
     \label{trefoilrecomb}
\end{center}
\end{figure}

As an example, view Figure~\ref{trefoilrecomb}. Here we see that two reconnections suffice to undo the trefoil knot. Each reconnection is interpreted as an oriented saddle cobordism.
Writhe is preserved throughout. The surface produced by the two reconnections has genus equal to one. Thus any single reconnection for the trefoil will have to produce a non-trivial link.
\bigbreak

The main mathematical result of this paper is the following Theorem, proved in Section 5.\\

\noindent {\bf Reconnection Theorem.} If $K$ is a positive unsplittable link diagram with with $c(K)$ crossings and $s(K)$ Seifert circles, then the reconnection number of $K$ is  given by the formula $$R(K) = c(K) - s(K) +1.$$\\

With this result we give an exact formula for the reconnection numbers for positive links and knots and thus give a strong basis for further experiments with physical vortices and mathematically modeled vortices.\\

\bigbreak

\section{Spanning Surfaces for Knots and LInks}

\begin{figure}
     \begin{center}
     \begin{tabular}{c}
     \includegraphics[width=7cm]{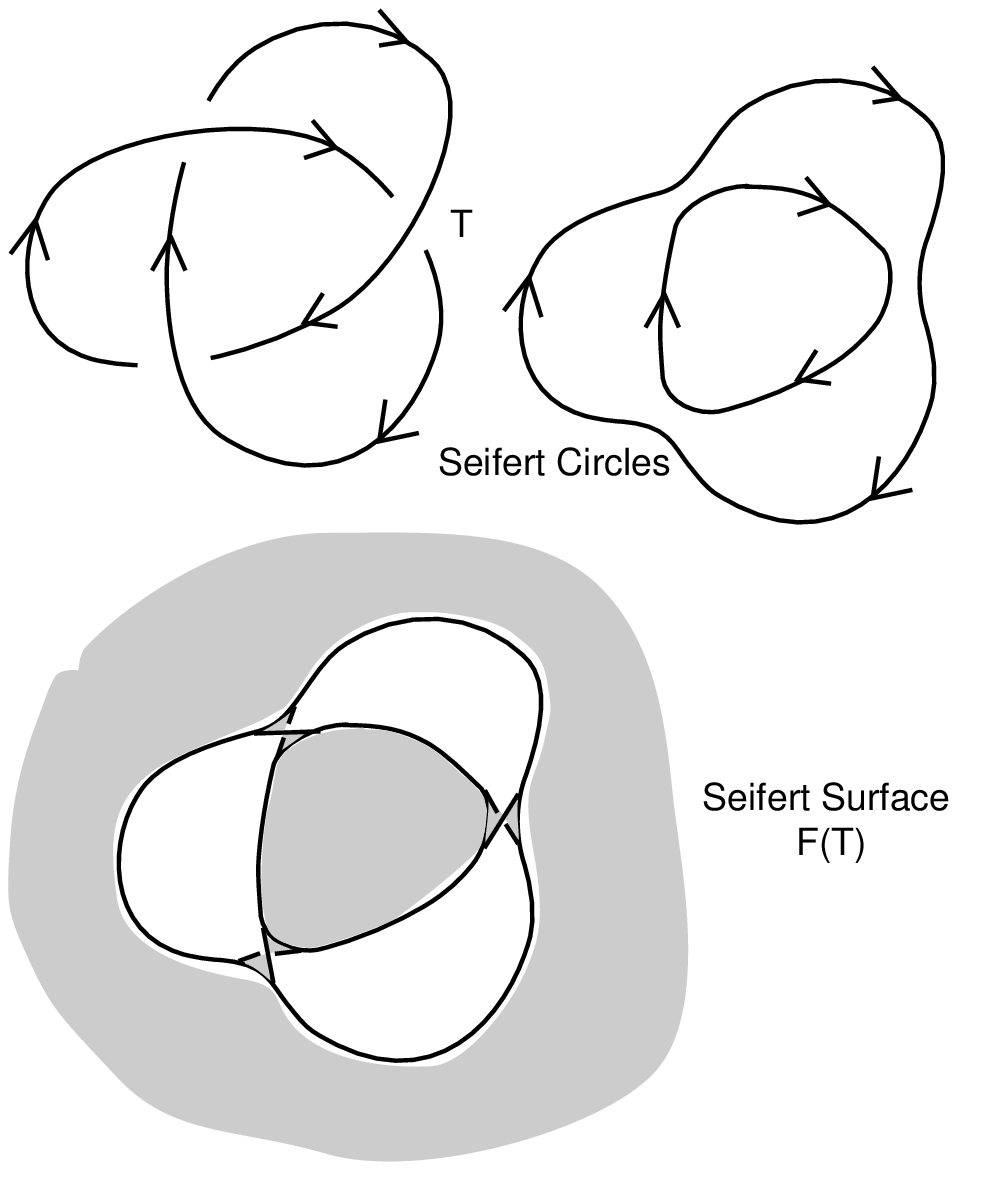}
     \end{tabular}
     \caption{\bf Classical Seifert Surface}
     \label{seifert}
\end{center}
\end{figure}

\begin{figure}
     \begin{center}
     \begin{tabular}{c}
     \includegraphics[width=7cm]{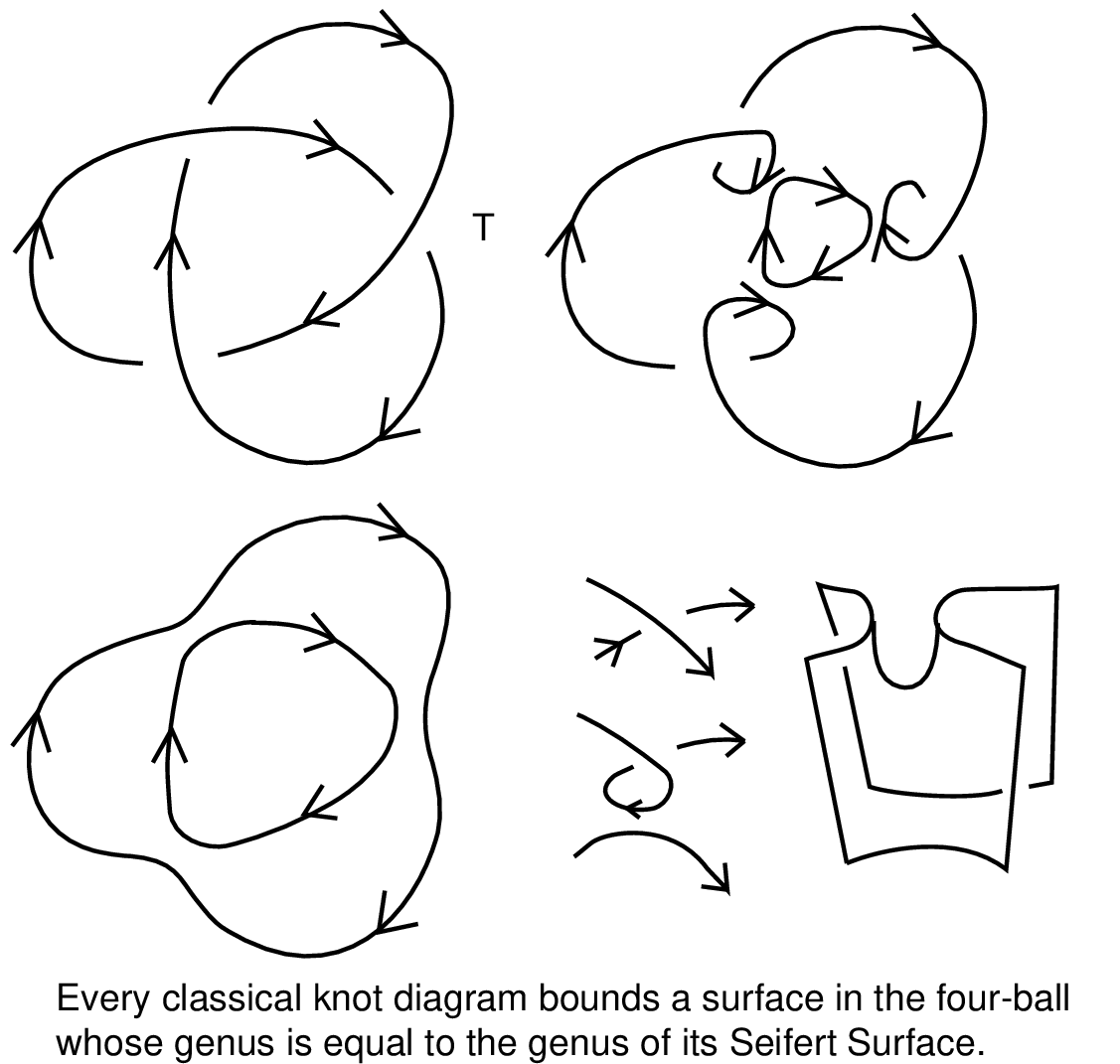}
     \end{tabular}
     \caption{\bf Classical Cobordism Surface}
     \label{classicalcob}
\end{center}
\end{figure}

\begin{figure}
     \begin{center}
     \begin{tabular}{c}
     \includegraphics[width=8cm]{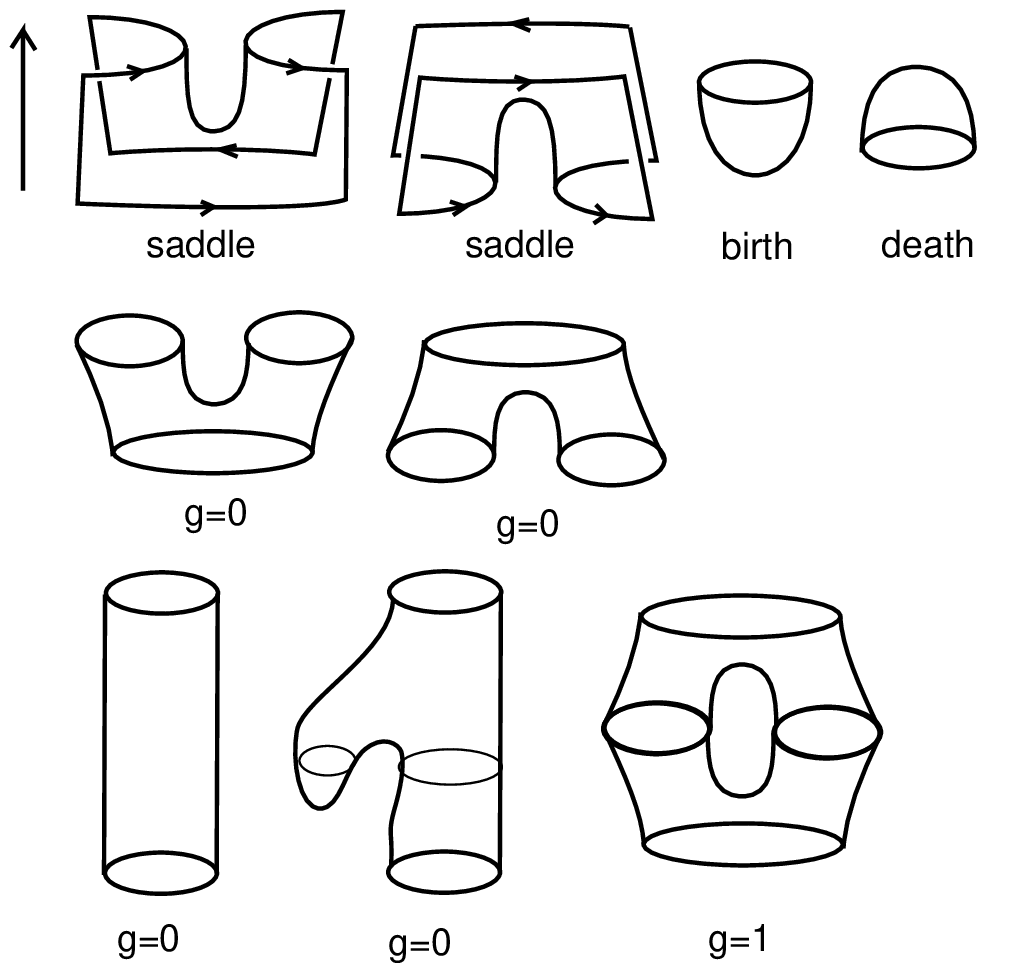}
     \end{tabular}
     \caption{\bf Saddles, Births and Deaths}
     \label{saddle}
\end{center}
\end{figure}

It is a well-known that every oriented classical (a one dimensional closed curve or curves embedded in three dimensional space)  knot or link bounds an embedded orientable surface in three-space. A representative surface of this kind can be obtained by the algorithm due to Seifert (See \cite{OK,FKT,KP}). We have illustrated Seifert's algorithm for a trefoil diagram in Figure~\ref{seifert}. The algorithm proceeds as follows: At each oriented crossing in a given diagram 
$K,$ smooth that crossing in the oriented manner (reconnecting the arcs locally so that the crossing disappears and the connections respect the orientation). The result of this operation is a collection of oriented simple closed curves in the plane, usually called the {\it Seifert circles}. To form the {\it Seifert surface} $F(K)$ for the diagram $K,$ attach disjoint discs to each of the Seifert circles, and connect these discs to one another by local half-twisted bands at the sites of the smoothing of the diagram. This process is indicated in the Figure~\ref{seifert}. In that figure we have not completed the illustration of the outer disc.
\bigbreak

It is important to observe that we can calculate the genus of the resulting surface quite easily from the combinatorics of the link diagram $K$ and knowledge of its number of components $\mu(K).$
\bigbreak

\noindent{\bf Lemma.} Let $K$ be a classical connected link diagram with $c(K)$ crossings and $s(K)$ Seifert circles.
then the genus of the Seifert Surface $F(K)$ is given by the formula
$$g(F(K)) =(1/2)( c(K) - s(K) + 1 - (\mu(K) -1))$$
where $c(K)$ denotes the number of crossings in the diagram $K,$ $s(K)$ denotes the number of Seifert circles for $K,$ and $\mu(K)$ denotes the number of components of the link $K.$
\bigbreak

\noindent {\bf Proof.} The surface $F(K),$ as described prior to the statement of the Lemma, is the homotopy type of a cell complex consisting of the projected graph of the knot diagram with 2-cells attached to each cycle in the graph that corresponds to a Seifert circle and 2-cells attached to each component of the link.  Thus we have that the Euler characteristic of this surface is given by the formula $$\chi(F(K)) = c - e + s + \mu$$ where $c,$ the number of crossings in the diagram, is the number of zero-cells, $e$ is the number of one-cells (edges) in the projected diagram (from node to node), $s$ is the number of Seifert circles and $\mu = \mu(K)$ as these are in correspondence with the 2-cells. We know that $4c = 2e$ since there are four edges locally incident to each crossing. Thus,
$$\chi(F(K)) =  - c + s + \mu$$ Furthermore, we have that $\chi(F(K)) = 2 - 2g(F(K)),$ since  this surface is orientable. From this it follows that $2-2g(F(K)) = -c + s + \mu,$ and hence
$$g(F(K)) =(1/2)(c - s+1 - (\mu -1)) =(1/2)( c(K) - s(K) + 1 - (\mu(K) -1)).$$  This completes the proof. \\
\bigbreak

We now observe that {\it for any classical link $K,$ there is a surface bounding that link in the four-ball that is homeomorphic to the Seifert surface}.  One can construct this surface by  pushing the Seifert
surface into the four-ball keeping it fixed along the boundary. We will give here a  different description of
this surface as indicated in Figure~\ref{classicalcob}. In that figure we {\it perform a reconnection at every crossing of the diagram.} The result is a collection of unknotted and unlinked curves. By our interpretation of surfaces in the four-ball obtained by saddle moves  (reconnections) and isotopies, we can then bound each of these curves by discs (via deaths of circles) and obtain a surface 
$S(K)$ embedded in the four-ball with boundary $K.$ As the reader can easily see, the curves produced by the saddle transformations are in one-to-one correspondence with the Seifert circles for
$K,$ and it easy to verify that $S(K)$ is homeomorphic with the Seifert surface $F(K)$ with boundary the link $K.$
Thus we know that $g(S(K)) =(1/2)( c(K) - s(K) + 1 -(\mu(K) -1)).$ Note that the genus of the surface is one half the rank of its first homology group,
\bigbreak

\section{The Seifert Pairing and the Signature}
Let $F$ be a spanning surface in three-space for a knot or link  $K.$ We define a linking number measure of the embedding of $F$ via the {\it Seifert paring} defined as an asymmetric bilinear form
$$\Theta: H_{1}(F) \times H_{1}(F) \longrightarrow \cal {Z},$$ given by the formula 
$$\Theta(x,y) = Lk(x^{*}, y)$$ where $Lk$ denotes linking number and $x^{*}$ denotes the result of translating the chain representing $x$ in the first homology group $H_{1}(F)$ along the positive normal to the surface $F$ so that $x^{*}$ is supported in the complement $S^{3} - K.$ The Seifert pairing gives invariant information about the topology of the embedding of the spanning surface for $K,$ and it also gives information about the embedding of $K.$  Below, we list properties of the Seifert pairing. See \cite{OK} for more information.\\

\begin{enumerate}
\item Let $\Delta_{K}(t) = Det(\Theta - t \Theta^{T})$ where $\Theta^{T} (x,y) = \Theta(y,x)$ is the {\it transpose} of $\Theta.$ Then $\Delta_{K}(t)$ is the {\it Alexander polynomial} of $K.$
\item Let $\sigma(K) = Signature (\Theta + \Theta^{T}).$ This is the {\it signature} of the link $K.$ It is the signature of the symmetric bilinear pairing
$$\Theta + \Theta^{T} :  H_{1}(F) \times H_{1}(F) \longrightarrow \cal {Z}.$$
A remarkable property of the signature of a link is that it forms a lower bound for the four-ball genus of that link \cite{Mura,KT}. That is, we have the following inequality.
$$|\sigma(K)|   \le 2g_{4}(L) +  \mu(K) - 1$$ where $g_{4}(L)$ denotes the least genus among connected orientable surfaces in the four-ball $D^{4}$ that span the link $K$ in the three-sphere $S^3$, seen as the boundary of the four-ball.
Here $\mu(K)$ is the number of components of the link $K.$ Thus, in the case of a knot we have $\mu(K) = 1$ and so $$|\sigma(K)| \le 2g_{4}(L).$$
\end{enumerate}

\section{Applications to Vortex Degeneration}

If $K$ is a knot or link representing a vortex that can undergo oriented writhe-preserving reconnection, then each act of reconnection can be interpreted as a saddle-point cobordism. Thus, if the knot undergoes
$N$ reconnections to produce a collection of unlinked circles, then we can interpret the cobordism as producing a surface in the 4-ball of genus no more than $N/2.$ \\

The relevance of the process view of cobordisms of knots is that the saddle passage is exactly a writhe-preserving oriented reconnection. In studying vortices one is interested in the least number of such reconnections that will transform a knot or link to a collection of unknotted circles. We define $R(K)$ to be the least number of reconnections that will transform $K$ to a collection of unknotted circles. When $K$ is a link, with more than one component, then we are concerned with the reconnection number when $K$ is $\it unsplittable.$ A link $K$ is said to be unsplittable if there is no isotopy in three dimensional space  that can carry the link into disjoint non-empty parts that are in separate three-balls. Thus we are only interested in links that will require reconnection in order to become separated into
a collection of unknotted circles. In speaking of an unsplittable link, we can include the case of knots, since no knot can be split apart into non-empty parts in disjoint three-balls.\\

Using the four-ball genus, we have:\\

\noindent {\bf Theorem.} Let $g_{4}(K)$ denote the least genus of four-ball spanning surfaces of an unsplittable link $K$ in the three-sphere. Let $R(K)$ be the reconnection number for $K.$ Then 
 $$2g_{4}(K) + \mu(K) -1 \le R(K).$$  For the signature of the link $K,$ we conclude that 
  $$|\sigma(K)|  \le R(K).$$ \\ 
 
 \noindent {\bf Proof.}  A generic connected surface that forms a cobordism between an unsplittable link $K$
 and an unknotted loop will require at least $\mu(K) -1$ saddle points (reconnections). Pairs of saddle moves contribute to the increment of one handle, whence an addition of one to the genus.
 This implies that twice the four-ball genus plus $\mu(K)-1$  is 
 a lower bound on the reconnection number.  Thus $2g_{4}(K) + \mu(K) -1 \le R(K).$
 The result about the signature follows from the fact, quoted in the last section, that 
$$|\sigma(K)|  \le 2g_{4}(L) +  \mu(K) - 1  $$
 where $\mu(K)$ is the number of components of the link $K.$\\

\noindent {\bf Remark.} The key point about this Theorem is that if we know the four-ball genus of a link $K,$ then we have a strong bound on the reconnection number $R(K)$ of $K.$
A case in point is the result of Rasmussen \cite{Rasmussen} (an application of Khovanov Homology. See also \cite{LKKho,DKK,Kho}.) that tells us the following:\\

\noindent{\bf Definition.} An oriented knot $K$  is said to be {\it positive} if all its crossings are of positive type. For reference, the crossings in the trefoil knot in Figure~\ref{classicalcob} are positive.\\

\noindent{\bf Theorem (Rasmussen \cite{Rasmussen}).} Let $K$ be a positive knot, then the four-ball genus of $K$ is equal to the genus of the Seifert surface for $K.$ Thus
$$ g_{4}(K) = \frac{1}{2}(c(K) - s(K) +1).$$\\

From this we obtain the following results about reconnection:\\

\noindent{\bf Theorem.} If $K$ is a positive knot, then the reconnection number of $K$ is bounded below by twice the genus of the Seifert spanning surface of $K.$ That is,
$$R(K) \ge 2g_{4}(K) = c(K) - s(K) +1$$ where $c(K)$ is the number of crossings of $K$ and $s(K)$ is equal to the number of Seifert circles of $K.$ \\

\noindent {\bf Proof.} This result follows from Rasmussen's result \cite{Rasmussen} that $g_{4}(K)$ is equal to the genus of a Seifert surface for $K$ and our discussion prior to the Theorem.\\

\noindent {\bf Theorem.} If $K$ is a positive knot with a positive diagram with $c(K)$ crossings and $s(K)$ Seifert circles, then the reconnection number of $K$ is exactly twice the genus of the Seifert spanning 
surface of $K.$ That is, we have the formula $$R(K) = c(K) - s(K) +1.$$\\

\noindent {\bf Proof.} Since we know that $R(K) \ge c(K) - s(K) +1,$ it suffices to show how to perform $c(K) - s(K) +1$ reconnections and obtain an unknotted loop. From our description of 
the construction of the Seifert surface we showed that one reconnection at each crossing in $K$ produces the set of Seifert circles. 
We can choose to reconnect at all but $s(K) -1$ crossings so that the result is a single unknotted loop. In order to understand the geometry of this reassembly it is helpful to phrase it more generally.
The set of Seifert circles is a special case of an arbitrary collection of disjoint circles in the plane. Given two circles, let a reassembly between them be denoted by a straight line segment drawn from one to
the other. Given $s$ circles in the plane, we can ask how many reassembly lines need to be drawn to create a simple closed curve. The answer, of course, is $s-1.$ In the general case, we are free to reassemble wherever we please, rather than at specific sites. It is easy to see that for the case of the Seifert circles for any link diagram the reassemblies can be chosen so that the resulting single loop is unknotted. This means that by performing $c(K) - (s(K) -1)$ reconnections on any connected link diagram, we can produce an unknotted curve. Such a reconnection is sufficient, for positive knots, to show that $R(K) = c(K) - s(K) +1.$
This completes the proof.\\

\noindent{\bf Remark.} The above result extends to positive links by using extensions of Rasmussen's work, as we shall see below. The work of Rasmussen for knots gives a combinatorial proof of  Milnor's conjecture \cite{Milnor} that the four-ball genus of a torus knot of type $(a,b)$ is $(a-1)(b-1)/2.$ Milnor's conjecture was first proved by Kronheimer and Mrowka \cite{KM1} by using
gauge theory. The construction in the proof of the Theorem above shows that, {\it for any knot or link $K,$ $R(K) \le c(K) - s(K) + 1.$ } In the case of positive links, this upper bound is also a lower bound.\\

\begin{figure}
     \begin{center}
     \begin{tabular}{c}
     \includegraphics[width=8cm]{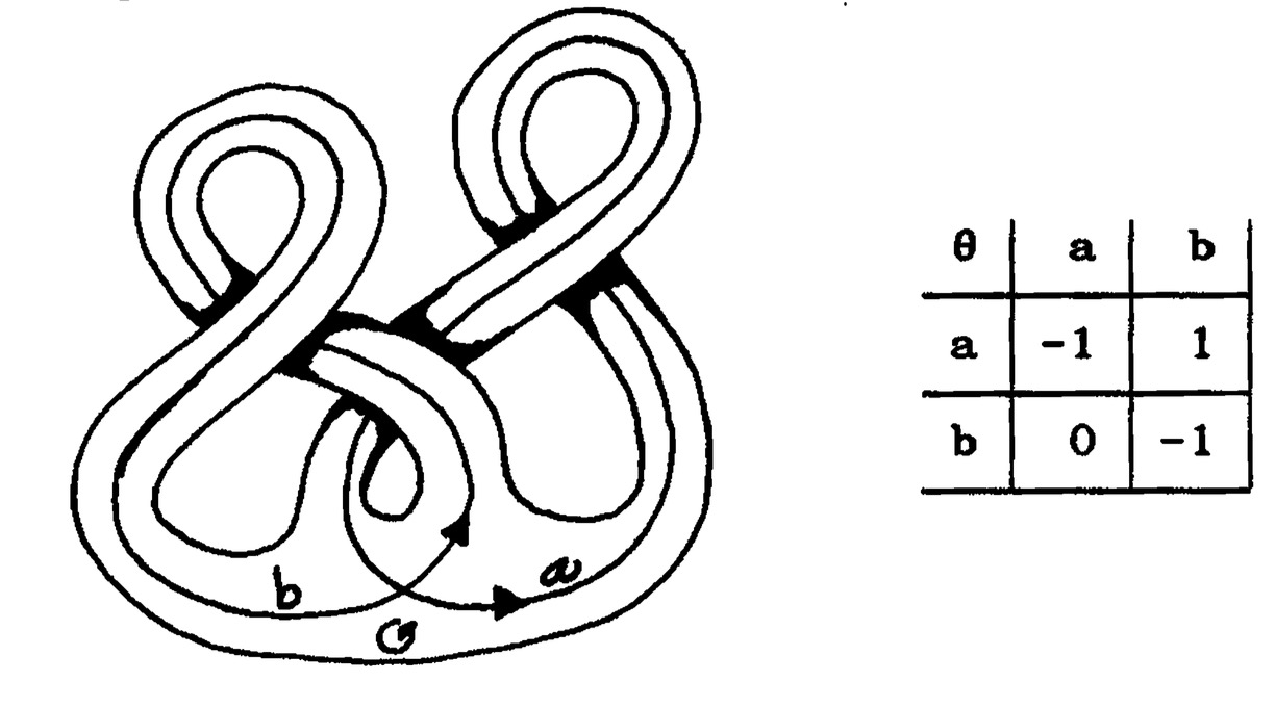}
     \end{tabular}
     \caption{\bf Seifert Pairing for Surface Bounding Trefoil Knot}
     \label{seifertpairing}
\end{center}
\end{figure}

\begin{figure}
     \begin{center}
     \begin{tabular}{c}
     \includegraphics[width=8cm]{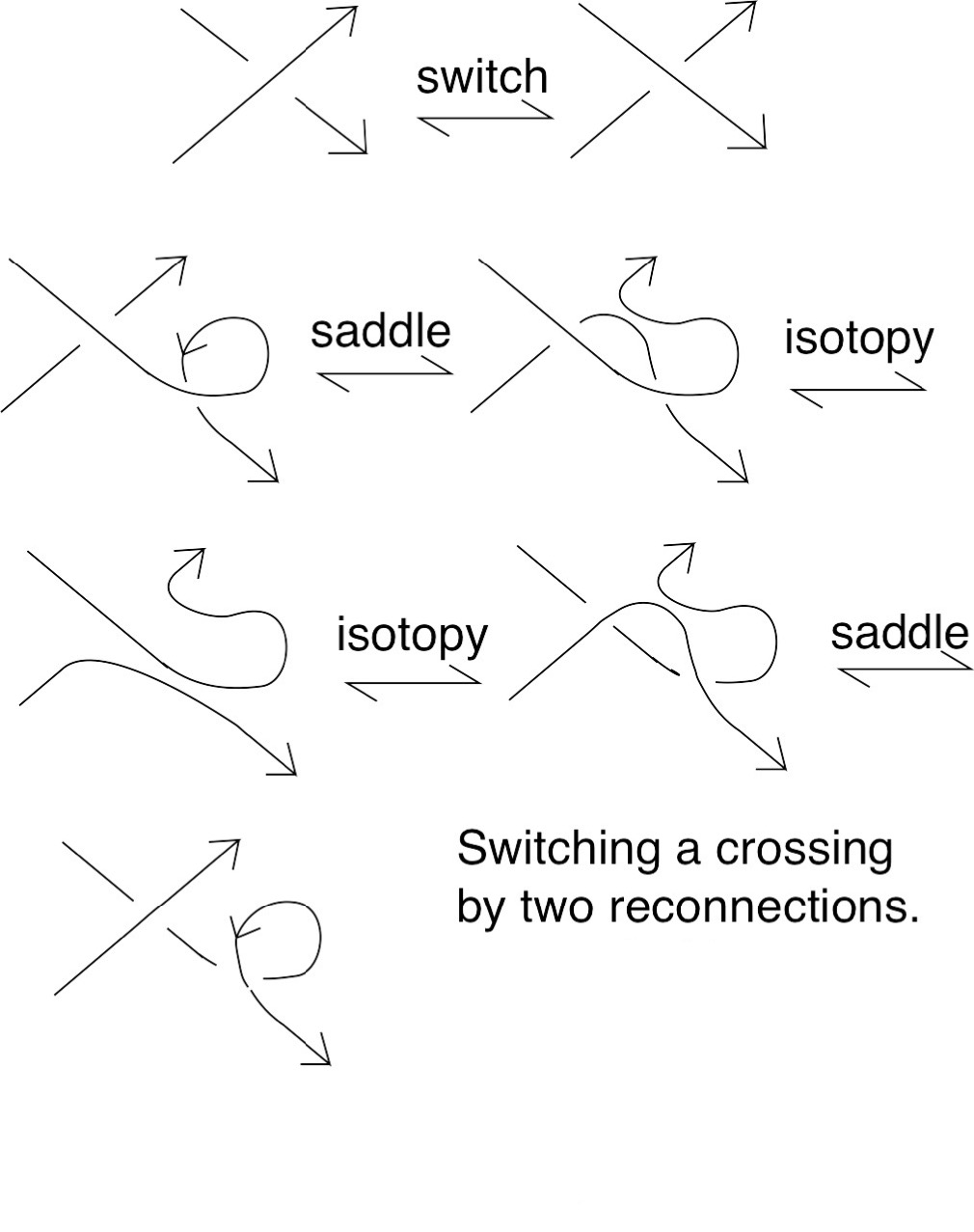}
     \end{tabular}
     \caption{\bf Crossing Switch in Two Reconnections}
     \label{crossing}
\end{center}
\end{figure}

\begin{figure}
     \begin{center}
     \begin{tabular}{c}
     \includegraphics[width=8cm]{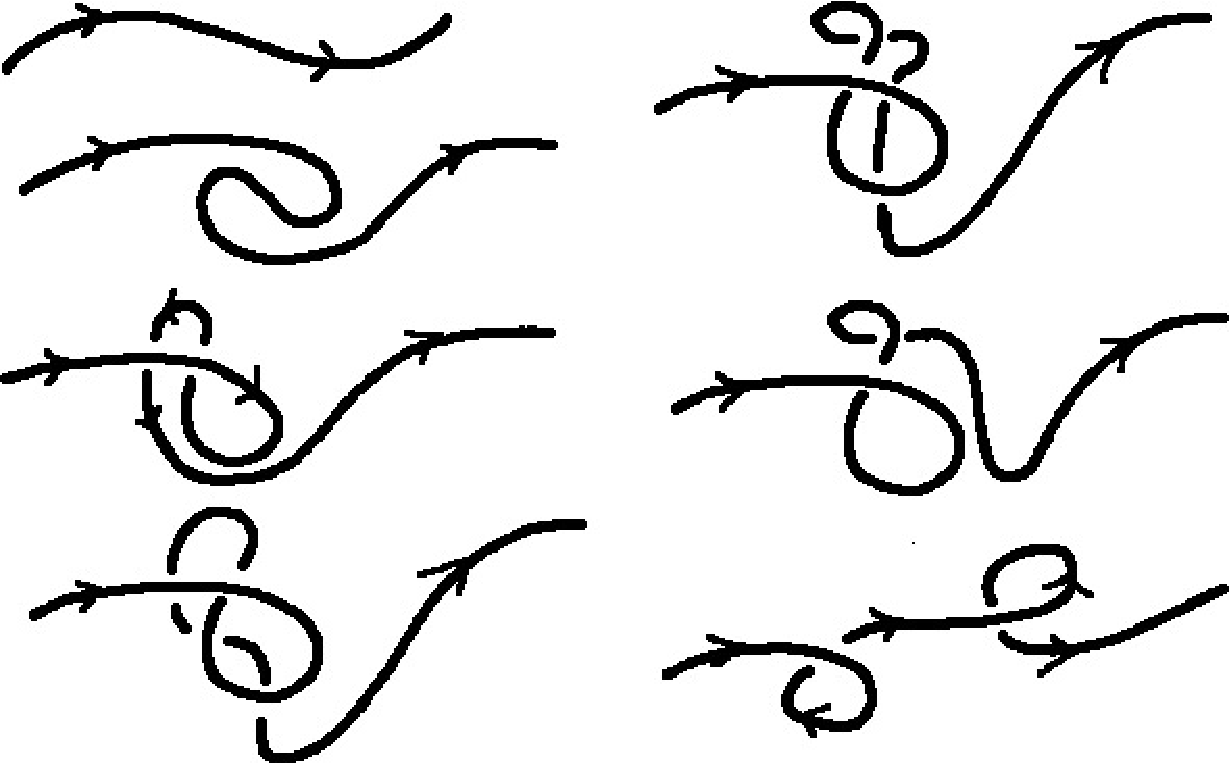}
     \end{tabular}
     \caption{\bf Producing Opposite Curls}
     \label{curls}
\end{center}
\end{figure}

\begin{figure}
     \begin{center}
     \begin{tabular}{c}
     \includegraphics[width=8cm]{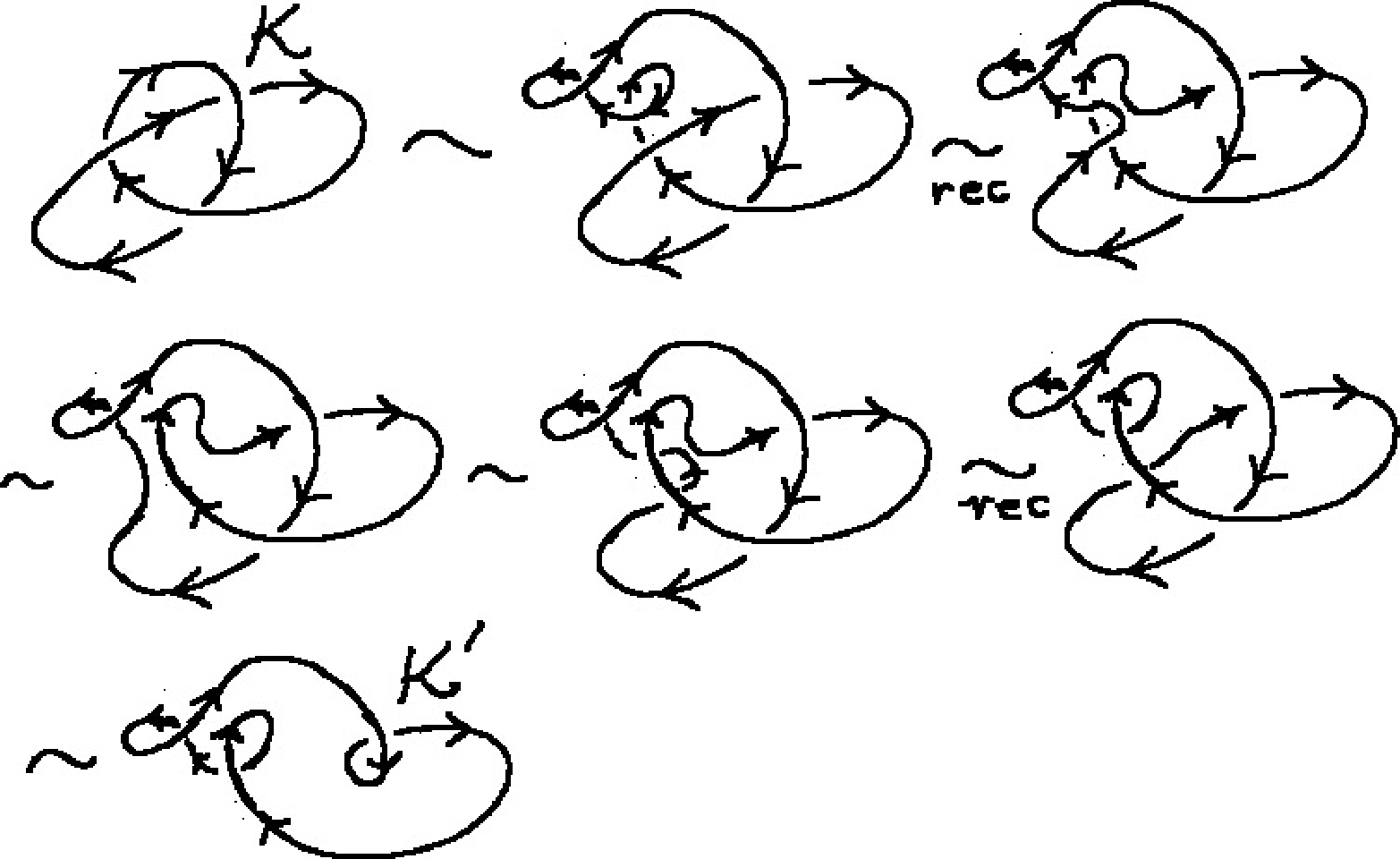}
     \end{tabular}
     \caption{\bf Trefoil Reconnection by Crossing Switch}
     \label{srecon}
\end{center}
\end{figure}

\begin{figure}
     \begin{center}
     \begin{tabular}{c}
     \includegraphics[width=10cm]{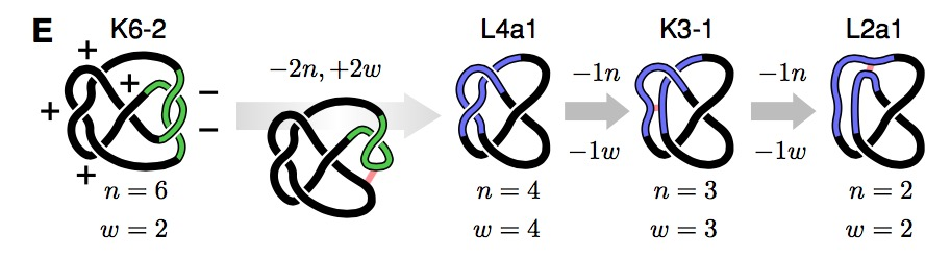}
     \end{tabular}
     \caption{\bf Reconnections for $6_2.$}
     \label{sixtwo}
\end{center}
\end{figure}

\begin{figure}
     \begin{center}
     \begin{tabular}{c}
     \includegraphics[width=8cm]{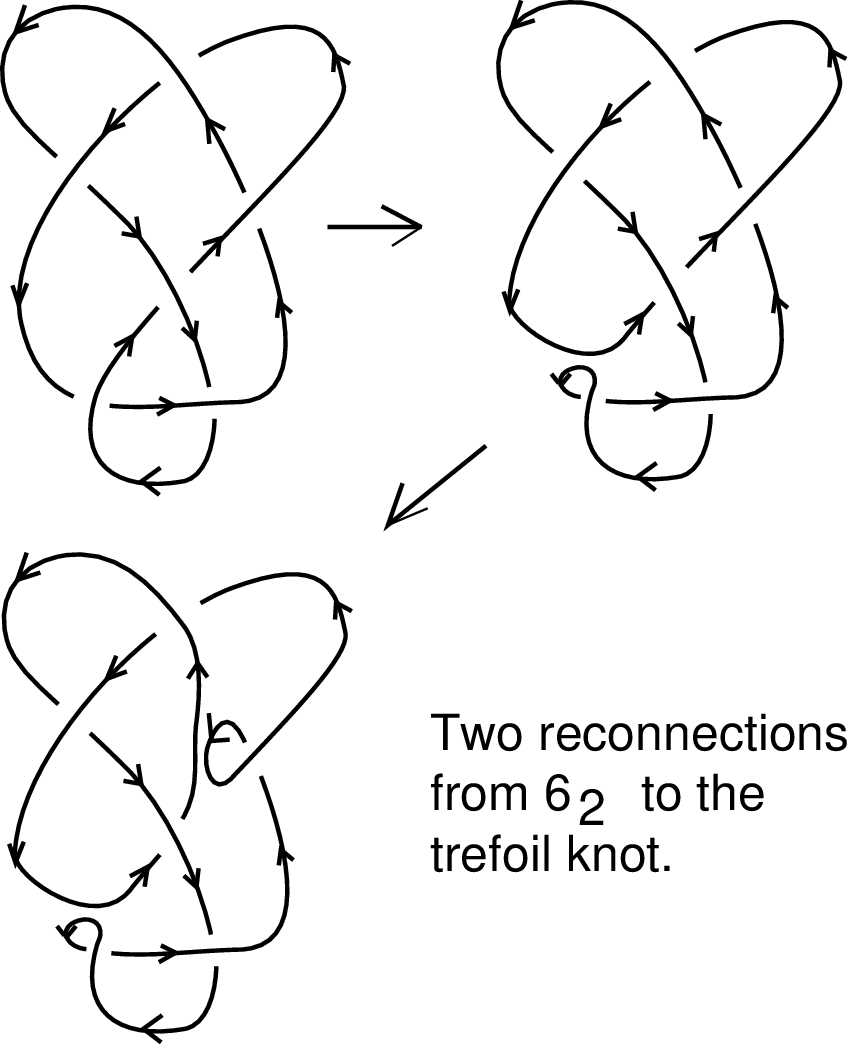}
     \end{tabular}
     \caption{\bf $6_2$ has two reconnections to the Trefoil. Hence four reconnections to an unknot.}
     \label{sixtwo1}
\end{center}
\end{figure}

\begin{figure}
     \begin{center}
     \begin{tabular}{c}
     \includegraphics[width=8cm]{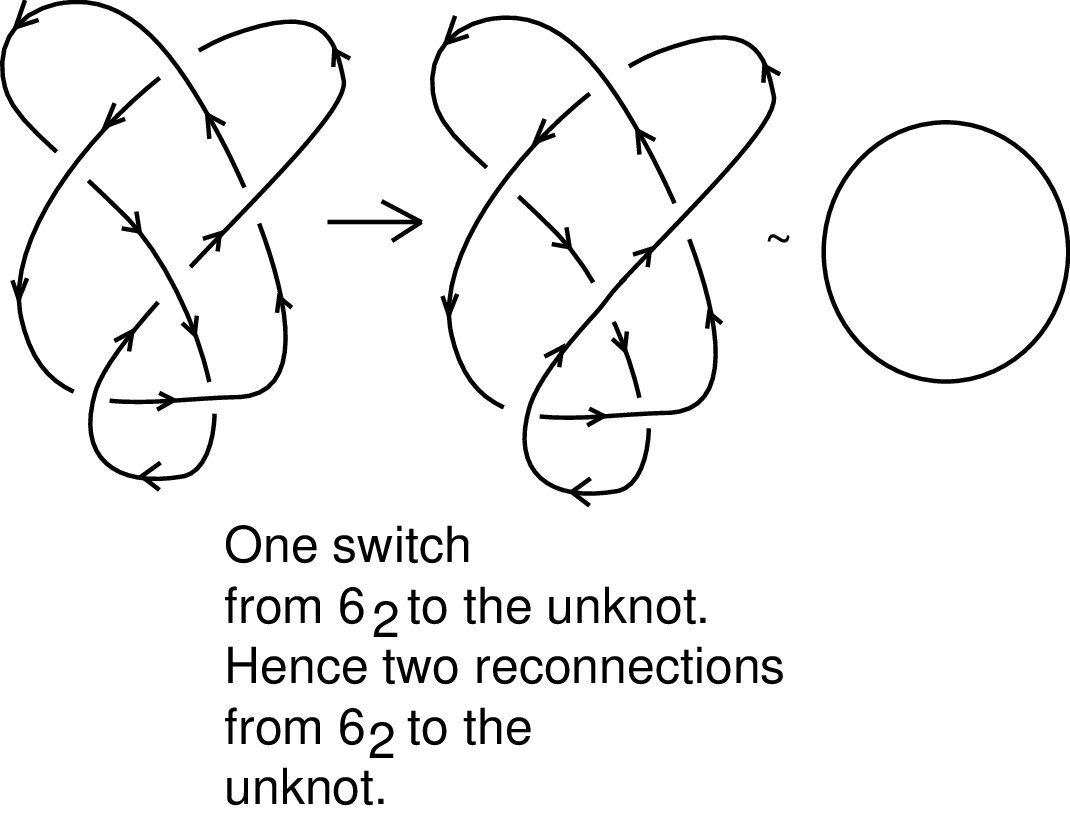}
     \end{tabular}
     \caption{\bf By Switch $6_2$ has two reconnections to the Unknot.}
     \label{sixtwoswitch}
\end{center}
\end{figure}

\begin{figure}
     \begin{center}
     \begin{tabular}{c}
     \includegraphics[width=8cm]{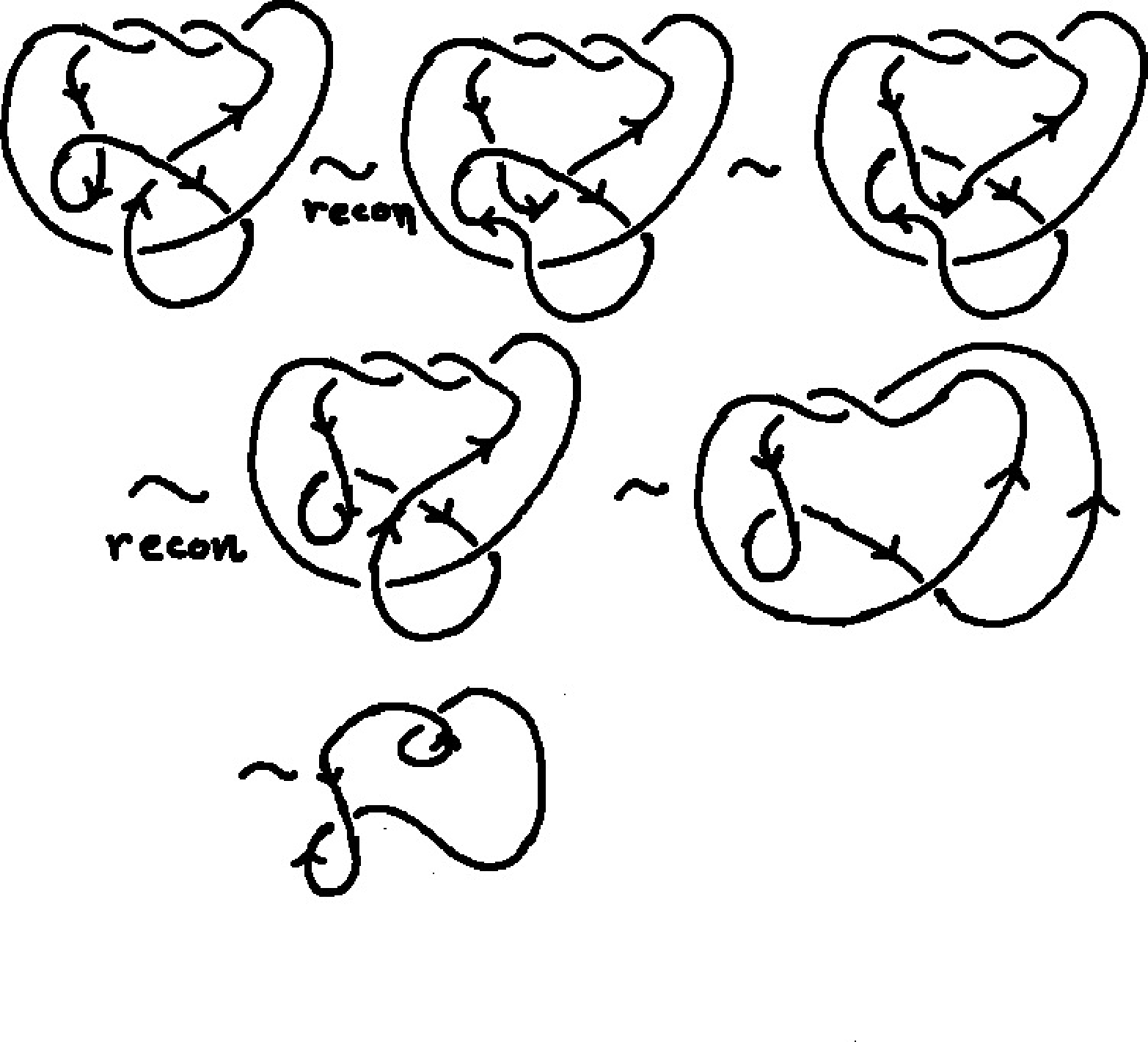}
     \end{tabular}
     \caption{\bf A two-fold reconnection for $6_2$}
     \label{sixtworecon}
\end{center}
\end{figure}

\begin{figure}
     \begin{center}
     \begin{tabular}{c}
     \includegraphics[width=10cm]{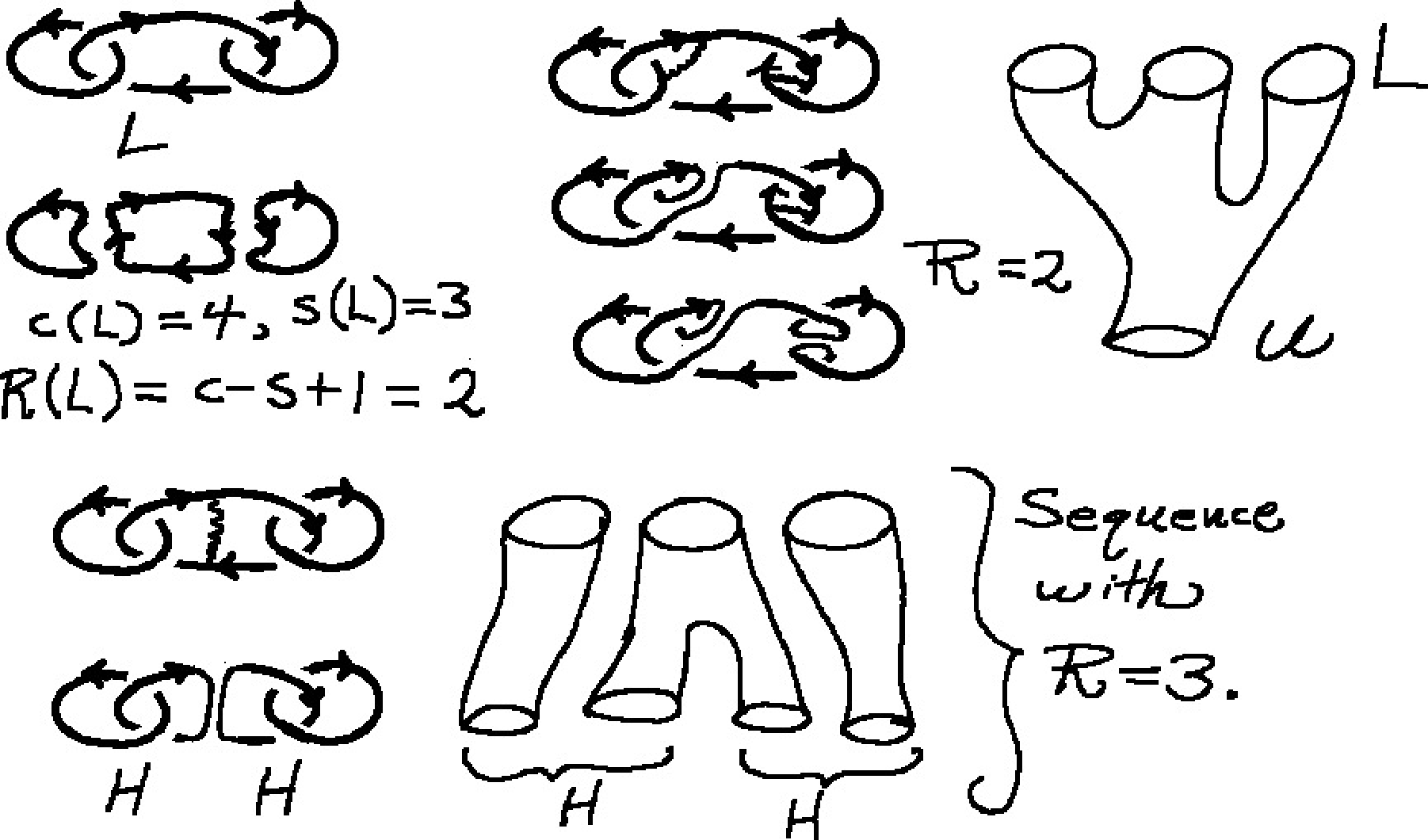}
     \end{tabular}
     \caption{\bf A Link Reconnection}
     \label{chain}
\end{center}
\end{figure}

\begin{figure}
     \begin{center}
     \begin{tabular}{c}
     \includegraphics[width=8cm]{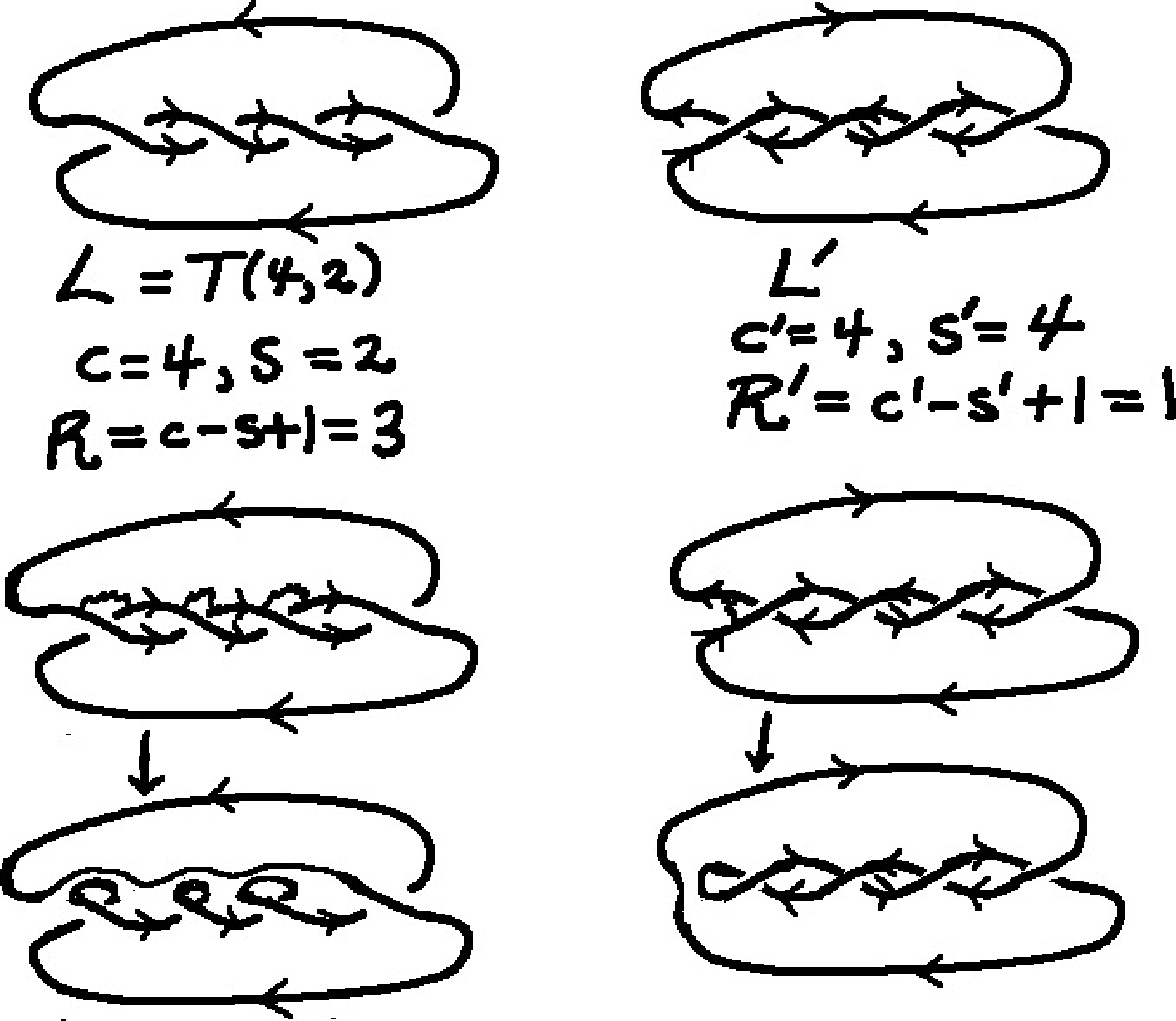}
     \end{tabular}
     \caption{\bf Reconnection dependent on orientation.}
     \label{orient}
\end{center}
\end{figure}

\begin{figure}
     \begin{center}
     \begin{tabular}{c}
     \includegraphics[width=10cm]{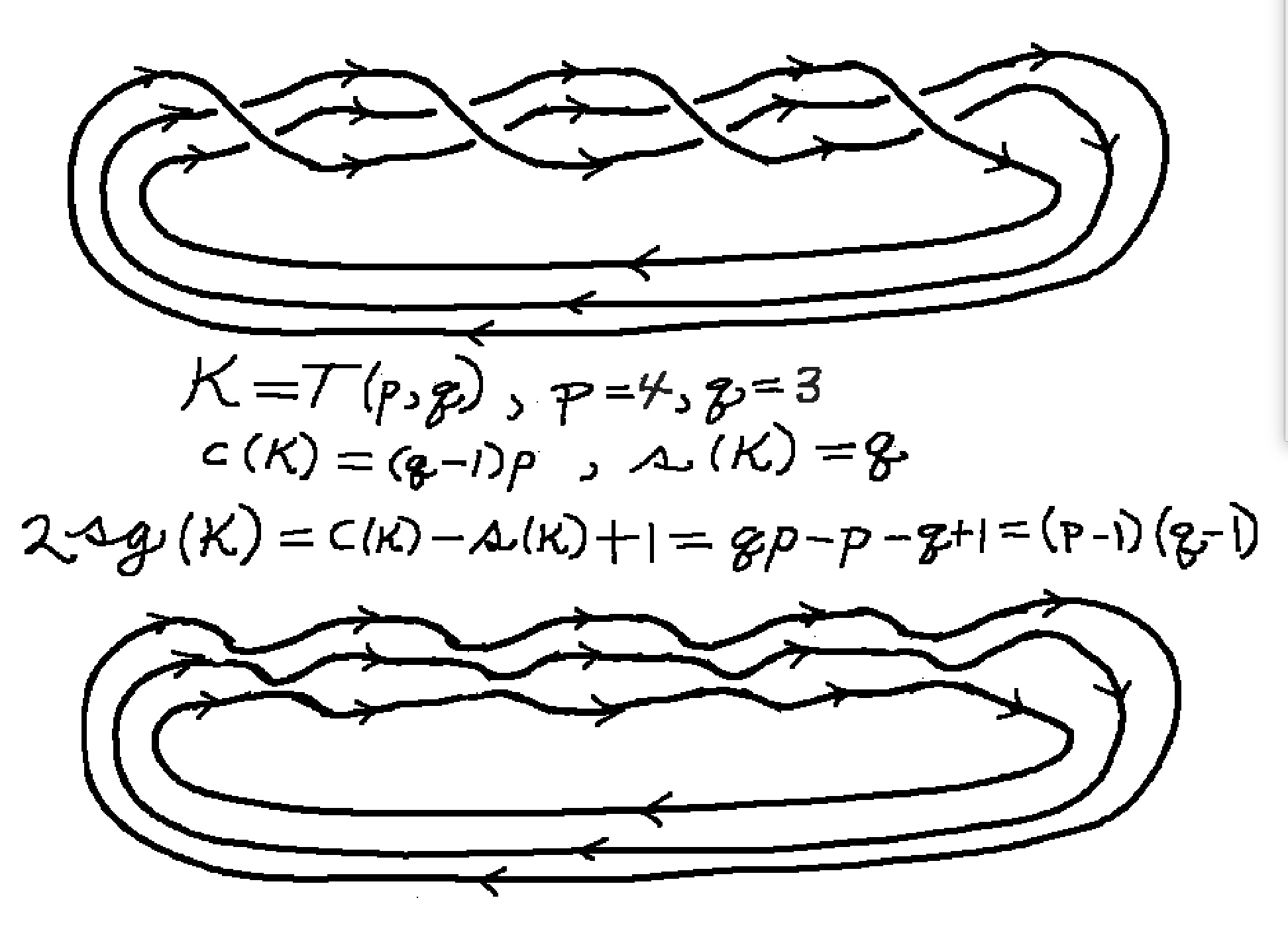}
     \end{tabular}
     \caption{\bf Seifert Genus of Torus Knots}
     \label{torus}
\end{center}
\end{figure}

\begin{figure}
     \begin{center}
     \begin{tabular}{c}
     \includegraphics[width=10cm]{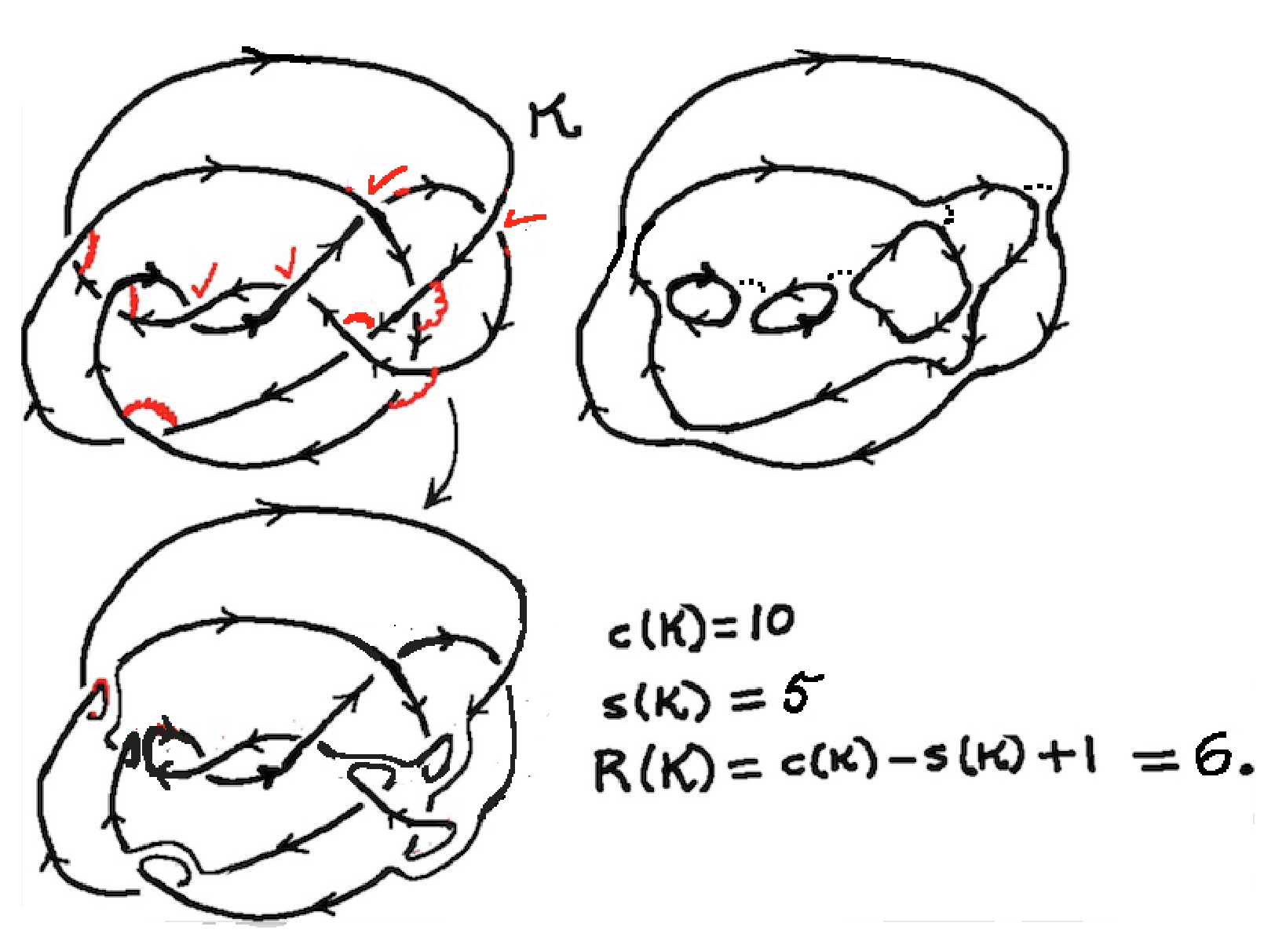}
     \end{tabular}
     \caption{\bf A Positive Knot}
     \label{positive}
\end{center}
\end{figure}

\begin{figure}
     \begin{center}
     \begin{tabular}{c}
     \includegraphics[width=10cm]{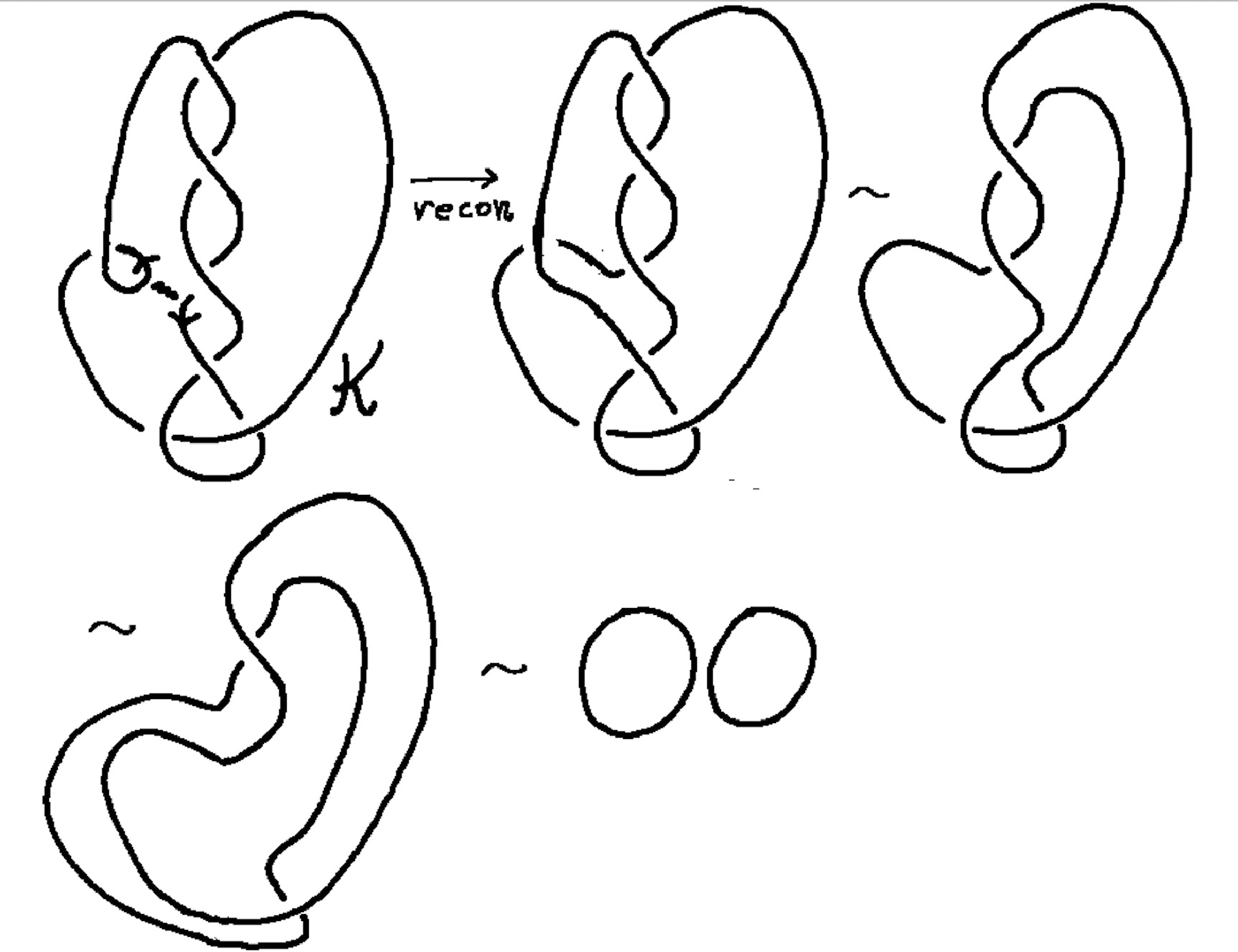}
     \end{tabular}
     \caption{\bf One reconnection suffices.}
     \label{onesaddle}
\end{center}
\end{figure}

\noindent {\bf Generalization to Links.} We now generalize the above results to the case of unsplittable links.\\

\noindent {\bf Reconnection Theorem.} If $K$ is a positive unsplittable link diagram with with $c(K)$ crossings and $s(K)$ Seifert circles, then the reconnection number of $K$ is  given by the formula $$R(K) = c(K) - s(K) +1.$$
Note that this formula generalizes our result for positive knots.\\

\noindent {\bf Proof.}  The result of Rasmussen for knots generalizes to links. This can be seen directly by examining the paper by Rasmussen, and it has been formalized by Cavallo \cite{Cavallo}.
See also \cite{BW,Abe}.
We state the result and refer the reader to Cavallo's paper for the details.
This is the result: For a positive unsplittable link, $ g_{4}(K) = (c(K) - s(K) + 1 - (\mu(K) -1))/2.$ The reader will note that this is the genus of the Seifert spanning surface for $K.$ Thus this statement is a direct generalization of the result of Rasmussen.  Thus we have that $$2g_{4}(K) + \mu(K) -1 = c(K) - s(K) +1,$$ and we know from our previous Theorem that 
$$2g_{4}(K) + \mu(K) -1 \le R(K).$$ Therefore $c(K) - s(K) + 1 \le R(K),$ and by the same constructive method we have already explained, one can perform $c(K) - s(K) + 1$ reconnections and transform $K$ to an unknotted loop.
Therefore $c(K) - s(K) + 1 = R(K).$ This completes the proof of the Theorem.\\

\noindent {\bf Example 1.} If $T$ is the trefoil knot, then $|\sigma(K)| = 2$, from which we conclude that it will take at least two reconnections to undo the trefoil, and indeed two reconnections suffice.
The Seifert Pairing for the trefoil knot is illustrated in Figure~\ref{seifertpairing}.
The matrix for the pairing is
$$\Theta = 	\left(\begin{array}{cc}
			-1&1\\
			 0& -1
			\end{array}\right).$$ 
Thus 
$$\Theta + \Theta^{T} =  \left(\begin{array}{cc}
			-2&1\\
			 1 &-2
			\end{array}\right),$$
from which it is easy to see that $|\sigma(K)| = 2.$\\

\noindent {\bf Example 2.}	In this example, we first point out that the switching of a crossing in a knot or link can be accomplished with two saddle moves that create a single addition to the genus in a cobordism 
surface. See Figure~\ref{crossing} for a detailed illustration of the process that gives this result. In this figure we show how to accomplish this switch by using a curl nearby the crossing. After the process both the crossing and the curl have been switched and the total writhe remains the same as at the beginning. In a knot where there is no nearby curl we can, in principle, obtain one for this process by producing two opposite curls from nothing by using only the second and third Reidemeister moves, as shown in Figure~\ref{curls}. 
This example shows that if we are looking for the reconnection number of a knot or link we can examine how many crossings it
will require to undo the knot. For example in the case of a trefoil knot, it requires one crossing switch to undo the trefoil and so we see that it can be undone by two reconnections as shown in Figure~\ref{srecon}.
Of course this method of undoing the trefoil by crossing switching is much more complex than our first method shown in Figure~\ref{recon}. Nevertheless, we can conclude that \\

\noindent {\bf Theorem} If a knot of link $K$ has unknotting number $u(K),$ then $R(K) \le 2u(K).$\\

Here the unknotting number $u(K)$ is the least number of crossing switches needed to unknot the knot or link $K.$ The unknotting number is an invariant of $K,$ but in general it is very difficult to compute this
number.\\

Here is another specific example of this type. Examine Figure~\ref{sixtwo}, Figure~\ref{sixtwo1} and Figure~\ref{sixtworecon}. In these figures we see that two natural reconnections take the knot $K = 6_2$ to a trefoil knot. Then two more reconnections will undo the trefoil knot. Thus we see that $K=6_2$ can be undone in $4$ reconnections. This is how it did happen using the Gross-Pitaevskii simulation in \cite{KKI}. However, $6_2$ can be unknotted with one crossing switch and so we conclude that $R(6_2) = 2$ even though the physical simulation suggested that it might be 4. Figure~\ref{sixtworecon} illustrates graphically how these two reconnections could occur in practice. Certainly we need to assess the probability of such sequences occurring in physical modes.\\

\noindent{\bf Example 3.} In Figure~\ref{chain} we show a positive link $L$ of three components with $c(L) = 4, S(L) = 3.$ Thus by our Reconnection Theorem we have that $R(L) = c(L) - s(L) + 1 = 4-3+ 1 = 2.$
The figure illustrates a reconnection sequence to an unknot with three reconnections as predicted by the Theorem. Note that in this case, the Seifert surface has zero genus. In general we have, for positive links, that 
$R(L) = g_{4}(L) + \mu(L) -1.$ In this case $g_{4}(L) = 0$ and $R(L) = \mu(L) -1 = 3-1 = 2.$ We also illustrate how a different choice of initial reconnection leads to a cascade into two Hopf links and the need for three reconnections along this pathway.\\

\noindent{\bf Example 4.} In Figure|\ref{orient} we show a link $L$ (distinct from the previous example) that is a torus link of type $(4,2)$ and another link $L'$ that is, as unoriented links, the mirror image of the first link.
The links $L$ and $L'$ are both positive but with different orientations. We find that $c(L) = c(L') = 4$ but $s(L) =2$ while $s(L')= 4.$ As a consequence $R(L) = 3$ while $R(L') = 1.$ The reconnections are indicated.
This example shows how the reconnection number is sensitive to orientation for links.\\

\noindent {\bf Example 5.}	Consider a torus knot $K = T(p,q)$ of type $(p,q)$ as shown in Figure~\ref{torus}.  The figure illustrates a (4,3) torus knot, but the pattern is the same in the general case.
One closes $q$ copies of a braid of $p$ strands consisting of positive crossings and representing a cyclic permutation of order $p$ so that the first strand goes over all the other strands to the last place, and each 
of the other strands goes up one strand-place. All crossings in the generating braid are positive. From this description it is easy to see, as in the figure, that there are $q$ Seifert circles so that $c(K) = (q-1)p$ and
$s(K)=q.$ Thus, from the Reconnection Theorem above, we have that $R(K) = c(K) - s(K) + 1 = (q-1)p -q + 1 =  pq - p - q+1 = (p-1)(q-1).$ This example alone should become the basis for much experimental research since it is
quite possible to prepare torus knots as vortices either with computer models or via the Gross-Pitaevskii equations.\\

\noindent {\bf Example 6.} As we have pointed out in the Reconnection Theorem, we know the reconnection number for any positive link. In Figure~\ref{positive} we illustrate one more example.
In this example, we have check-marked a set of crossings at which no reconnection is performed. There are $s(K)-1$ checked crossings. By performing reconnections at the remaining 
$c(K) - s(K) + 1$ crossings we obtain an unknot. This is an illustration of the method of the Reconnection Theorem that produces reconnection sequences and thereby proves that he reconnection number of any positive knot is equal to $c(K) - s(K) + 1.$  Since we now know the exact reconnection number of an infinite class of knots, including all the torus knots, this provides a useful ground for 
physical experimentation.\\

\noindent {\bf Example 7.} View Figure~\ref{onesaddle} where we show a knot $K$ (the stevedore's knot) that becomes an unlinked link by a single reconnection. The stevedore's knot is indeed knotted and has unknotting number equal to one. This is the simplest example we know showing that the reconnection number can be strictly less than twice the unknotting number. There are many examples of this phenomenon. \\

\noindent {\bf Remark.} Any given knot or link can be seen to give rise to various cascades of reconnection, leading to unknots and unlinks. Our formulas for the reconnection numbers of positive links give
a way to navigate these pathways. It should be remarked that in \cite{Ricca1,Ricca2} patterns of cascade in relation to Conway and Homflypt polynomials have been noted. Since it is the case \cite{FKT} that the genus of the Seifert spanning surface for positive knots and links is seen in degrees of the Conway polynomial, some of these patterns are related to our exact calculations of the reconnection numbers. Again, more experimental 
work is needed in this domain.\\

\noindent{\bf Remark.} The results in this section all generalize essentially verbatim to virtual knot theory via our results \cite{DKK} generalizing the Rasmussen Invariant for virtual knot theory.
We can define the reconnection number $R(K)$ of a virtual knot or link diagram to be the least number of oriented saddle moves that produces an unknot or unlink in the virtual category.
It follows from our results that the reconnection number for a positive virtual knot or link is given by the same formula as in the Reconnection Theorem in the classical case. Since virtual knot theory
is a way to study knots and links in thickened surfaces, it is possible that experiments can be done with knotted vortices that are restricted to such three dimensional ambient spaces.\\

\end{document}